\documentclass[11pt]{amsart}
\usepackage{amssymb, latexsym, amsmath}
\usepackage{graphicx}
\usepackage{color}
\usepackage{showlabels}

\usepackage{datetime}

\begin{document}
\title[Quasi-local energy]{Quasi-local energy in presence of gravitational radiation}
\author{Po-Ning Chen, Mu-Tao Wang, and Shing-Tung Yau}
\date{\today}
\begin{abstract}
We discuss our recent work \cite{Chen-Wang-Yau4} in which gravitational radiation was studied by evaluating the Wang-Yau quasi-local mass of surfaces of fixed size at the infinity of both axial and polar perturbations of the Schwarzschild spacetime, \`{a} la Chandrasekhar \cite{Chandrasekhar}. \end{abstract}

\thanks{P.-N. Chen is supported by NSF grant DMS-1308164, M.-T. Wang is supported by NSF grant  DMS-1405152,  and S.-T. Yau is supported by NSF
grants  PHY-0714648 and DMS-1308244. This work was partially supported by a grant from the Simons Foundation (\#305519 to Mu-Tao Wang). Part of this work was carried out
when P.-N. Chen and M.-T. Wang were visiting the Department of Mathematics and the Center of Mathematical Sciences and Applications at Harvard  University.}

\address{Department of Mathematics\\
Columbia University\\ New York\\ NY 10027\\ USA}
\email{mtwang@math.columbia.edu}

\date{\usdate{\today}}

\maketitle

We compute the Wang-Yau quasi-local mass \cite{Wang-Yau1, Wang-Yau2} of ``spheres of unit size" at null infinity to capture the information of gravitational radiation.  The set-up, following Chandrasekhar \cite{Chandrasekhar}, is a gravitational perturbation of the Schwarzschild solution, which is governed by the Regge-Wheeler equation (see below).  We take a sphere of a fixed areal radius and push it all the way to null infinity. The limit of the geometric data is that of a standard configuration and thus the optimal embedding equation \cite{Wang-Yau1, Wang-Yau2, Chen-Wang-Yau1} can be solved.

  Let us first  consider the axial perturbations. The metric perturbation is of the form:
\[-(1-\frac{2m}{r}) dt^2+\frac{1}{1-\frac{2m}{r}} dr^2+r^2 d \theta^2+r^2 \sin^2\theta (d\phi-q_2 dr-q_3 d\theta)^2.\] The linearized vacuum Einstein equation is solved by a separation of variable Ansatz in which $q_2$ and $q_3$ are explicitly given by the Teukolsky function and the Legendre function. 

In particular, \[q_3= \sin(\sigma t)  \frac{C_\mu(\theta)}{\sin\theta} \frac{(r^2-2mr)}{\sigma^2 r^4} \frac{d}{dr} (rZ^{(-)})\] for a solution of frequency $\sigma$ and a separation of variable constant $\mu$.
Here $C_\mu(\theta)$ is related to the $\mu$-th Legendre function $P_\mu$ by 
\[
C_\mu(\theta) = \sin \theta \frac{d}{d\theta} (\frac{1}{\sin \theta} \frac{d P_\mu(\cos \theta)}{d \theta}).
 \]
After the change of variable \[r_*=r+2m\ln (\frac{r}{2m}-1),\] $Z^{(-)}$ satisfies the Regge-Wheeler equation:
\[(\frac{d^2}{dr_*^2}+\sigma^2)Z^{(-)}=V^{(-)}Z^{(-)},\] where
\[V^{(-)}=\frac{r^2-2mr}{r^5}[(\mu^2+2)r-6m],\] and $\mu$ is a separation of variable constant. 

On the Schwarzschild spacetime \[-(1-\frac{2m}{r}) dt^2+\frac{1}{1-\frac{2m}{r}} dr^2+r^2 d \theta^2+r^2 \sin^2\theta d\phi^2,\] we consider an asymptotically flat 
Cartesian coordinate system $(t, y_1, y_2, y_3)$ with $y_1=r\sin\theta\sin\phi, y_2=r\sin\theta\cos\phi, y_3=r\cos\theta$. Given $(d_1, d_2, d_3)\in \mathbb{R}^3$ with $d^2=\sum_{i=1}^3 d_i^2$ , consider  the 2-surface
\[\Sigma_{t,d}=\{(t, y_1, y_2, y_3) :\sum_{i=1}^3 (y_i-d_i)^2=1\}.\] We compute the quasi-local mass of $\Sigma_{t,d}$ as $d\rightarrow \infty$.

\begin{center}
\includegraphics[width=4in]{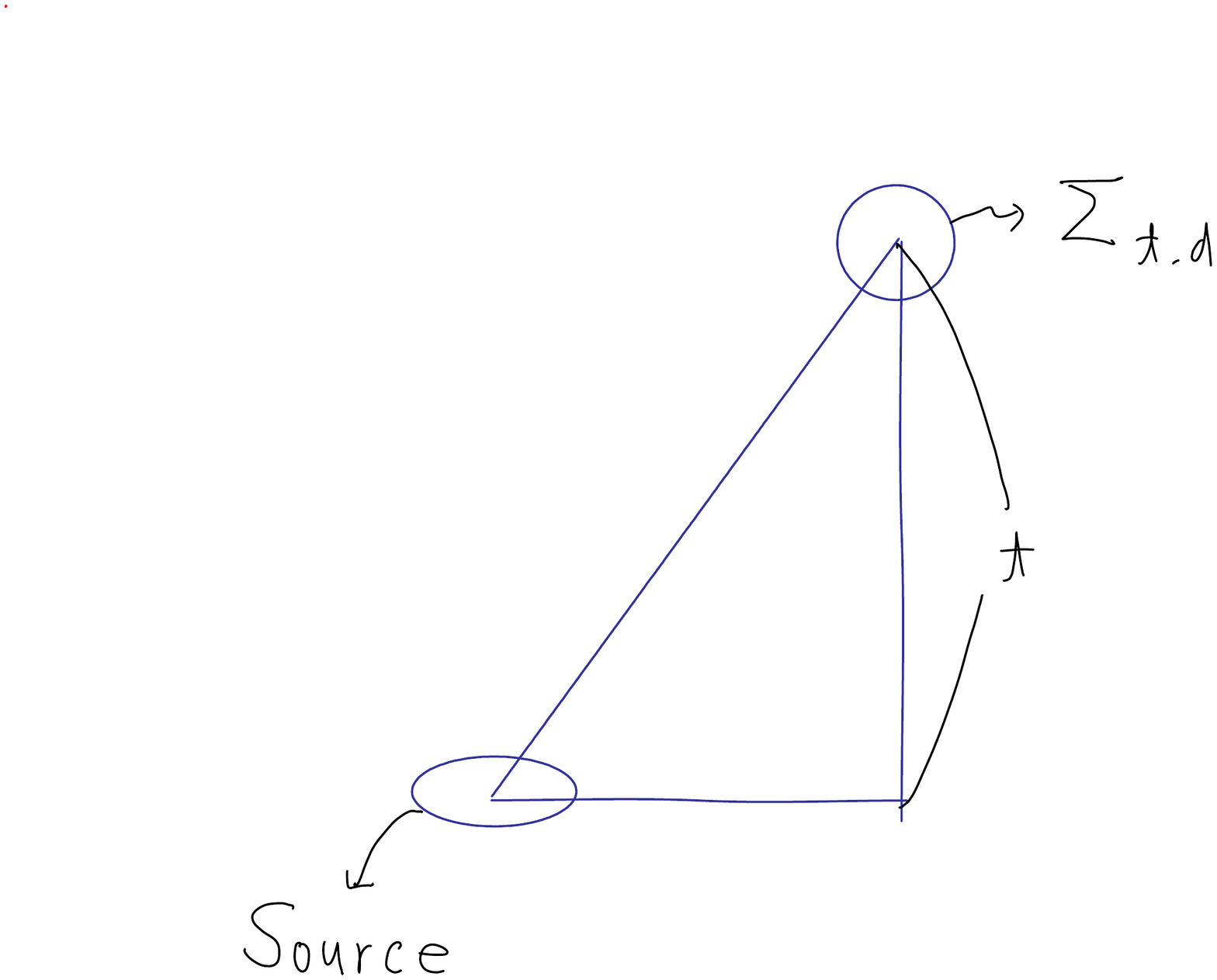}

\end{center}

Denote \[A({r})=\frac{(r^2-2mr)}{\sigma^2 r^3} \frac{d}{dr} (rZ^{(-)}).\] The linearized optimal embedding equation of $\Sigma_{t,d}$ is reduced to two linear elliptic equations on the unit
2-sphere $S^2$:
\[\begin{split} \Delta(\Delta+2)\tau&=[-A''(1-Z_1^2) +6A' Z_1+12A]Z_2 Z_3\\
(\Delta+2) N&=(A''-2A'Z_1+4A)Z_2Z_3, \end{split}\]
where $\tau$ and $N$ are the respective time and radial components of the solution, and $Z_1, Z_2, Z_3$ are the three standard first eigenfunctions of $S^2$. $A'$ and $A''$ are derivatives with respect to $r$, and $r^2$ is substituted by $r^2=d^2+2Z_1+1$ in the above equations.

The quasi-local mass of $\Sigma_{t,d}$ with respect to the optimal isometric embedding is then 
\[\begin{split}&E(\Sigma_{t,d}) =  C^2\{\sin^2(\sigma t )E_1+\sigma^2\cos^2(\sigma t)E_2\}+O(\frac{1}{d^3}),
\end{split}\] where $E_1$ and $E_2$ are two integrals on the standard unit 2-sphere, that depend on the solution $\tau$ and $N$ of the optimal isometric embedding
equation. Explicitly,
\[\begin{split} E_1&=\int_{S^2} (1/2)\left[A^2Z_2^2(7Z_3^2+1)+2AA'Z_1Z_3^2(3Z_2^2-1)-N(\Delta+2)N\right]\\
E_2&=\int_{S^2} \left[A^2 Z_2^2Z_3^2-\tau\Delta(\Delta+2)\tau\right].\end{split} \] In particular,
\[  \partial_t E(\Sigma_{t,d}) =  \frac{\sigma \sin ( 2 \sigma t )C^2(\theta)}{d^2}\{E_1- \sigma^2 E_2\}+O(\frac{1}{d^3}).  \] 
Let us compare the quasi-local mass on the small spheres $\Sigma_{t,d}$ along a certain direction to the quasi-local mass of the large coordinate spheres $S_{t,r}$.

\begin{center}
\includegraphics[width=4in]{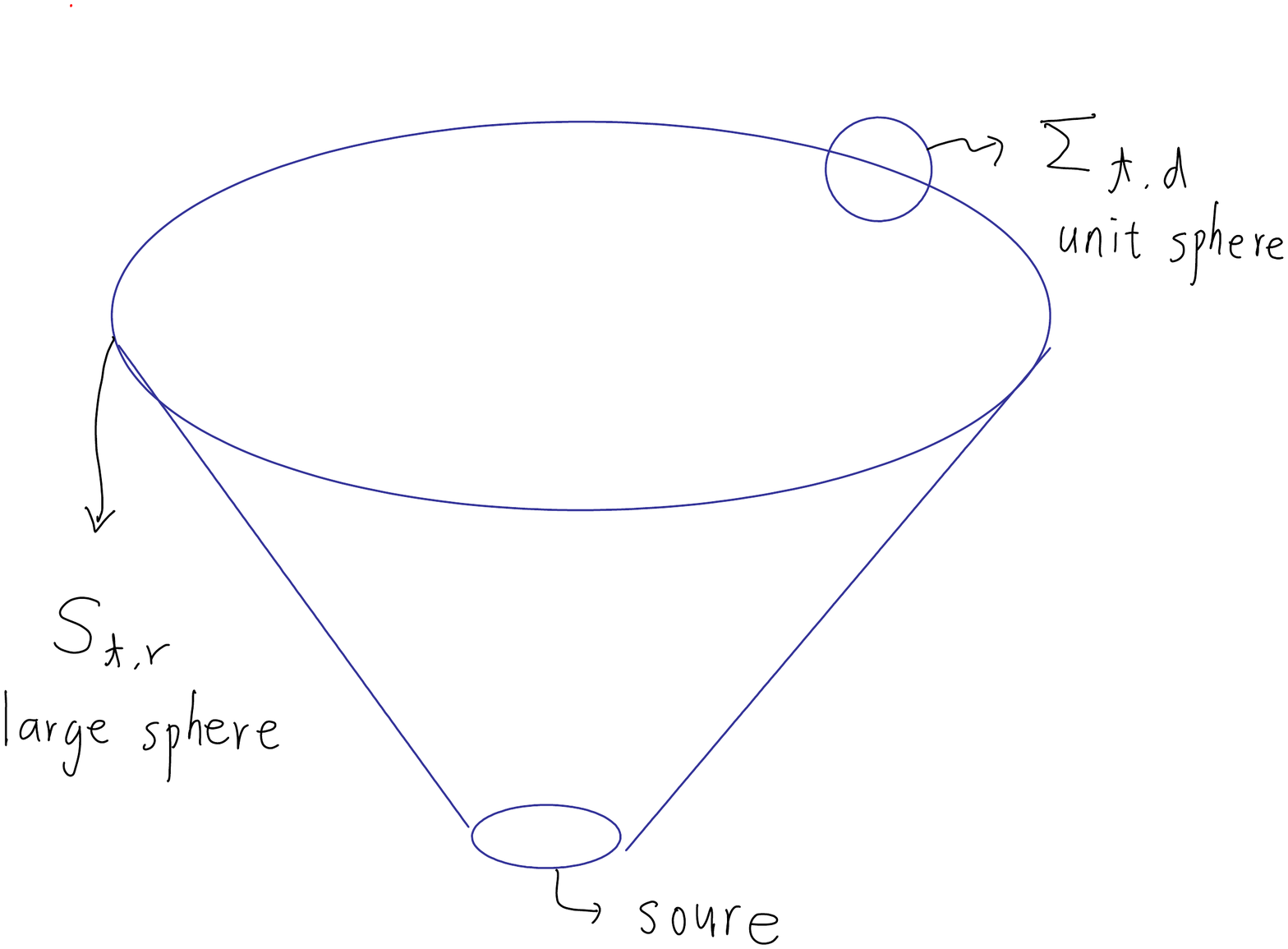}

\end{center}

Naively, one may expect to recover  $\partial_t E(S_{t,r})$ by integrating the energy radiated away at all directions $\partial_t E(\Sigma_{t,d}) $. However, our calculation indicates that  there are nonlinear correction terms from the quasi-local energy that should be taken into account.

We can also consider the polar perturbation of the Schwarzschild spacetime in which the metric coefficients $g_{tt}$, $g_{rr}$, $g_{\theta\theta}$, and $g_{\phi \phi}$ are perturbed in 
\[-(1-\frac{2m}{r}) dt^2+\frac{1}{1-\frac{2m}{r}} dr^2+r^2 d \theta^2+r^2 \sin^2\theta d\phi^2.\]

The gravitational perturbation is governed by the Zerilli equation
\[(\frac{d^2}{dr_*^2}+\sigma^2)Z^{(+)}=V^{(+)}Z^{(+)},\] where
\[V^{(+)}=\frac{2(r^2-2mr)}{r^5(nr+3m)^2}[n^2(n+1) r^3+3m n^2r^2+9m^2 nr+9m^3],\] and $n$ is the separation of variable constant.  Again, we compute the quasi-local mass of spheres of unit-size at null infinity. The calculation is similar to the axial perturbation case but the result is  different as the leading term is of the order $\frac{1}{d}$ (as opposed to $\frac{1}{d^2}$ for axial-perturbation) with nonzero coefficients. If such a linear perturbation can be realized as an actual perturbation of the Schwarzschild spacetime, the result would contradict the positivity of the quasi-local mass \cite{Liu-Yau, Wang-Yau1, Wang-Yau2}. From this, we deduce the following conclusion: There does not exist any gravitational perturbation of the Schwarzschild spacetime that is of purely polar type
in the sense of Chandrasekhar \cite{Chandrasekhar}.

 For an actual gravitational perturbation of the Schwarzschild solution, the vanishing of the $\frac{1}{d}$ gives a limiting integrand that 
integrates to zero on the limiting 2-sphere at null infinity. In fact, the quasi-local mass density $\rho$ (see \cite[equation 2.2]{Chen-Wang-Yau3}) of $\Sigma_{t,d}$ can be computed at the pointwise level. Up to an $O(\frac{1}{d^3})$ term
\[\begin{split} \rho&=(K-\frac{1}{4} |H|^2)\\
&-\frac{(|H|-2)^2}{4}
+\frac{1}{d^2}\{\frac{1}{2}|\nabla^2 N|^2+((\Delta+2)N)^2-\frac{1}{4}(\Delta N)^2\\
&-\frac{1}{4}(\Delta\tau)^2+\frac{1}{2}[\nabla^a\nabla^b(\tau_a\tau_b)-|\nabla \tau|^2
-\Delta |\nabla\tau|^2]\},\end{split}\] where $K$ is the Gauss curvature  of $\Sigma_{t,d}$. The first line, which integrates to zero, is of the order of $\frac{1}{d}$ and is exactly the mass aspect function of the Hawking mass \cite{Hawking}.  The $\frac{1}{d^2}$ term of the quasi local mass $\int_{\Sigma_d} \rho\,\, d\mu_{\Sigma_{t,d}}$ has contributions from the second and third 
lines (of the order of $\frac{1}{d^2}$), the $\frac{1}{d^2}$ term of the first line, and the $\frac{1}{d}$ term of the area element $d\mu_{\Sigma_{t,d}}$. The above integral formula is obtained after
performing integrations by parts and applying the optimal embedding equation several times. 

To each closed loop on the limiting 2-sphere at null infinity, we can thus associate a non-vanishing arc integral
that is of the order of $\frac{1}{d}$, where $d$ is the distance from the source.  We expect the freedom in varying the shape of the loop can increase the detectability of gravitational waves.

\end{document}